\documentclass[12pt]{amsart}
\usepackage{amsmath,amssymb}
\usepackage{amscd}
\usepackage{mathrsfs}
\usepackage[all]{xy}
\usepackage{pstricks}
\UseComputerModernTips
\usepackage{graphicx}
\DeclareGraphicsExtensions{.eps} {\bf }

\theoremstyle{plain}
\newtheorem{theorem}{Theorem}[section]

\newtheorem{lemma}{Lemma}[section]

\newcommand{\edge}[1]{\ar@{-}[#1]}
\theoremstyle{definition}
\newtheorem{definition}{Definition}[section]

\theoremstyle{remark}
\newtheorem{remark}{Remark}[section]

\numberwithin{equation}{section} \numberwithin{equation}{section}
\textheight=8.75in \theoremstyle{plain}
\begin{document}
\title[Limits of zeros]{Limits of zeros of polynomial sequences}
\author{Xinyun Zhu and George Grossman}
\address{Department of Mathematics\\ Central Michigan University\\ Mount Pleasant, MI 48859}
\email{zhu1x@cmich.edu,\, gross1w@cmich.edu} \date{\today}
\begin{abstract}  In the present paper we consider
$F_k(x)=x^{k}-\sum_{t=0}^{k-1}x^t,$ the characteristic polynomial of
the $k$-th order Fibonacci sequence, the latter denoted $G(k,l).$ We
determine the limits of the real roots of certain odd and even
degree polynomials related to the derivatives and integrals of
$F_k(x),$ that form infinite sequences of polynomials, of increasing
degree. In particular, as $k \rightarrow \infty,$  the limiting
values of the zeros are determined, for both odd and even cases. It
is also shown, in both cases, that the convergence is monotone for
sufficiently large degree. We give an upper bound for the modulus of
the  complex zeros of the polynomials for each sequence. This gives
a general solution related to problems considered by Dubeau 1989,
1993, Miles 1960, Flores 1967, Miller 1971 and later by the second
author in the present paper, and Narayan 1997.
\end{abstract}
\maketitle
{
\section{Introduction} The current work arose from consideration of
sequences of polynomials  ~\cite{george2} related to the asymptotic
behavior of their zeros. It is based on the following infinite
sequence of polynomials denoted as $\{F_k(x)\}_{k=1}^{\infty}$ for
convenience in the present paper which for $k \ge 2$, comprise the
characteristic polynomials of the $k$-th order Fibonacci sequence,
denoted by $G(k,l)$ where for $l>k\ge 2,$
\begin{equation*}
G(k,l)=\sum_{t=1}^{k}G(k,l-t),
\end{equation*}and $G(k,1)=1, \ G(k,t)=2^{t-2}, \ t=2,3,\ldots,k.$ For $k=2$ we obtain the well-known
Fibonacci sequence,
$\{1,1,2,3,5,8,\ldots,F_{n-1}+F_{n-2}=F_n,\ldots\}$.

It is also well-known that
\begin{equation*} \lim_{k\rightarrow
\infty}\frac{G(k,l+1)}{G(k,l)}=\phi_k, \ k \ge 2,
\end{equation*}
where $\phi_k$ is the positive zero of $F_k.$ Number theoretic
results concerning $G(k,l)$ are in ~\cite{george1}. A fractal
described by A. Dias, in A. Posamentier and I. Lehman's new book
~\cite{truss} was first published in ~\cite{george1}. The
significance of this fractal with respect to the present paper is
that the fractal dimension is $\ln({\phi_2})/\ln{2}.$

Miles 1960, ~\cite{miles} showed that the zeros of the sequence of
polynomials $\{F_{k}(x)\}, k \ge 2$ are distinct, all but one lies
in the unit disk and the latter is real and lies in the interval
$(1,2)$. Miller ~\cite{miller}, 1971 gave a different, shorter proof
of this result. Flores 1967,  ~\cite{flores},  showed that
$\phi_{k}\rightarrow 2$ monotonically as $k \rightarrow +\infty$ as
did Dubeau,  ~\cite{fd1}, ~\cite{fd2}. In ~\cite{george2} the
sequences $\{F_{k}^{'}(x)\}$ and $\{F_{k}^{''}(x)\}$ were studied
and we reproduce the following table for understanding and
motivation:
\begin{table}[h]
\begin{center}
\caption{Does Interval Contain a Root, yes or no?} \vspace{6pt}
\begin{tabular}{|l|c|c|c|c|c|c|c|} \hline  int/fn
&$F_{2k}$&$F_{2k+1}$&$F_{2k}^{\prime}$&$F_{2k+1}^{\prime}$&$F_{2k}^{\prime
\prime }$&$F_{2k+1}^{\prime \prime }$ \\ \hline
$(-1,0)$&yes&no&no&yes&yes ($k>1$)&no \\ \hline $(0,1]$&no&no&yes
($k=1$)&yes ($k=1$)&yes ($k=2$)&yes ($k=1,2$) \\ \hline
$(1,2)$&yes&yes&yes ($k>1$)&yes ($k>1$)&yes ($k>2$)&yes ($k>2$) \\
\hline
\end{tabular} \end{center} \end{table}

For the particular particular cases we find that
$F_{3}¥^{\prime}(1)=0,F_{2}^{\prime}¥(1/2)=0, \ \ F_{2}^{\prime
\prime}¥ = 2,  \ F_{3}¥^{\prime \prime}(1/3)=0, \ \ F_{5}¥^{\prime
\prime}(1)=0,\ \ F_{4}¥^{\prime \prime}((1+\sqrt{11/3})/4)=0.$

Note that in table 1, the number of negative roots is either $0$ or
$1$  for odd and even degree respectively, while there is always a
positive root in $(1,2)$ (for sufficiently large degree.) It was
indicated in ~\cite{george2} as an open question as to whether this
happens for higher derivatives and conjectured in \cite{george3}.

In ~\cite{george2} it was also shown that
$\lim_{k\rightarrow\infty}\theta_{k}= -1$ where $\theta_k$
is the negative zero of each term in $\{F_{2k}\}, \ k \ge 1.$
Similarly, by examining approximations to zeros, the same asymptotic
result was shown to hold for the sequences $\{F_{k}^{'}(x)\}$ and
$\{F_{k}^{''}(x)\}.$

In ~\cite{george3} a conjecture was also made concerning the real
zeros of the of $l$-th derivatives of each member of the sequence
$\{F_k\}_{k=2}^\infty$. Namely, the zeros of
$\{F_k^{(l)}\}_{k=2}^\infty$ exhibit the same (monotonic) behavior.
A conjecture that the complex zeros are all within the unit circle
was also made.

In this paper the question in ~\cite{george2} is answered, as are
the first two questions of ~\cite{george3}, affirmatively. The cases
of the complex zeros is still open, although we obtain an upper
bound. The present work also answers the same questions and yields
similar results for the $l$-th integral of $\{F_k\}$.

In the present paper then, we consider the following sets of
infinite sequences of polynomials given by,
\begin{equation*} U=\{ \{F_1, F_2, \ldots,\}, \ \{F^\prime_1, F^\prime_2,
\ldots,\},\{F^{\prime\prime}_1, F^{\prime\prime}_2,
\ldots,\},\ldots, \},
\end{equation*} and,
\begin{equation*}V=\left\{ \{F_1, F_2, \ldots,\},  \left\{\int F_1 dx, \int F_2 dx,
\ldots,\right\} ,\left\{\int\int F_1 dxdx, \int\int F_2 dxdx,
\ldots,\right\},\ldots,\right\}\end{equation*} where $F_1(x)=x-1$ and
\begin{equation*}F_k(x)=x^{k}-\sum_{t=0}^{k-1}x^t,k \ge 2.\end{equation*}

The sets $U, \ V$ are related to certain recurrence relations
~\cite{george4}, ~\cite{george5} having solutions that lead to
combinatorial identities. These recurrence relations result from a
factorization of $F_k(x),$ with unknown coefficients. Several
combinatorial identities are in ~\cite{sumid1}, ~\cite{fib2},
~\cite{sumid2}, for example it is shown in  ~\cite{sumid2} that for
any $c \neq -1,0$
\begin{eqnarray}\label{1.2}
\frac{\frac{1}{c^{2(n+1)}}-1}{1+c}&=&
\frac{1}{c^{n+2}}\sum_{i=1}^{n+1}\binom{n+i}{2i-1}\frac{(1-c)^{2i-1}}{c^{i-1}}
\nonumber \\
 &=&-\frac{1}{c}+\frac{1}{c^2}+\dots+\frac{1}{c^{2(n+1)}}, \ n \geq
0.
\end{eqnarray}\\\vspace{3mm}

 If $c \rightarrow -1$ in (\ref{1.2}) one obtains,
\begin{eqnarray*}2(n+1)=\sum_{i=1}^{n+1}\binom{n+i}{2i-1}2^{2i-1}(-1)^{n+i+1},
\end{eqnarray*} which is equivalent to a result in G. P\'{o}lya and G.
Szeg\"{o}, ~\cite{polszeg}.

The outline of the paper is as follows: in the next sections,
\S\ref{U}, \S\ref{V}, we give the three main results with proofs supported
in several lemmas. The first result deals with the set of
derivatives $U.$ The first and second derivative cases were
treated in ~\cite{george2}; the second and third results deal with
the set of integrals $V.$ The second result deals with the first
integral for which the proof leads to the general case and so is
included for interest and clarity of exposition.
\section{Results}\label{res}

\subsection{$U$ or derivative case}\label{U}
Now we consider the infinite sequence of polynomials $\{F_k^{(l)}(x)\}$ of the $l$-th derivative of the  sequence   $\{F_k(x)\}.$
\begin{definition}
We specify the following degree $j$ polynomial $D_j(x)$ to correspond with the $l$-th derivative of $F_{j+l}(x).$
\begin{equation}\label{equationFk}
D_j(x)=F_{j+l}^{(l)}(x)=l!\left(\binom{j+l}{l}x^j-\sum_{t=0}^{j-1}\binom{t+l}{l}x^t\right), \quad j\ge 1,
\end{equation}
 with $D_{0}(x)=l!.$
\end{definition}
\begin{lemma}
The $l$-th derivative of $F_k(x)$
is given by,
\begin{equation}\label{equationDk-l}
D_{k-l}(x)=\frac{\sum_{t=0}^{l+1}(-1)^{t}a_tx^{k+1-t}+(-1)^{l}l!}{(x-1)^{l+1}},\quad x\not= 1,
\end{equation}
 \noindent where each $a_i$ is a degree $l$ polynomial in $k$ with positive leading coefficient.
\end{lemma}
\begin{proof}
We can write
\begin{equation*}
F_k(x)=\frac{x^{k+1}-2x^{k}+1}{x-1}.
\end{equation*}
 We obtain the first derivative of $F_k(x)$ given by
\begin{equation*}
F'_k(x)=\frac{kx^{k+1}-(3k-1)x^k+2kx^{k-1}-1}{(x-1)^2},\quad x\not= 1.
\end{equation*}
\noindent Hence the statement is true for $l=1.$ Suppose the
statement is true for  $1\le l\le j.$  We have
\begin{equation} \label{equationDkj}
D_{k-j}(x)=\frac{\sum_{t=0}^{j+1}(-1)^{t}a_tx^{k+1-t}+(-1)^{j}j!}{(x-1)^{j+1}},
\end{equation}
\noindent where each $a_i$ is a degree $j$ polynomial in $k$ with
positive leading coefficient.

We obtain  the next derivative of (\ref{equationDkj}):
\begin{eqnarray}\label{equationDkjl}
D_{k-j-1}(x)
&=&\frac{\left(\sum_{t=0}^{j+1}(-1)^{t}a_t(k+1-t)x^{k-t}\right)(x-1)^{j+1}}{(x-1)^{2j+2}}\\
&-&\frac{(j+1)(x-1)^{j}\left(\sum_{t=0}^{j+1}(-1)^{t}a_tx^{k+1-t}+(-1)^{j}j!\right)}{(x-1)^{2j+2}}\nonumber\\
&=&\frac{\sum_{t=0}^{j+2}(-1)^{t}b_tx^{k+1-t}+(-1)^{j+1}(j+1)!}{(x-1)^{j+2}},\nonumber
\end{eqnarray}\\
where
\begin{eqnarray*}
b_{0}=a_{0}(k+1)-a_{0}(j+1)=a_0(k-j).
\end{eqnarray*}
For $1\le t\le j+1,$ we obtain by comparing the coefficients of like powers of $x$ in  (\ref{equationDkjl})
\begin{eqnarray*}
b_t&=&a_{t}(k+1-t)+a_{t-1}(k+2-t)-(j+1)a_t
\\ &=&a_t(k-j-t)-a_{t-1}(k+2-t),
\end{eqnarray*}
and
\begin{equation*}
b_{j+2}=a_{j+1}(k-j).
\end{equation*}
Hence the lemma follows.
\end{proof}
\begin{lemma}\label{droot} If $k-l$ is odd, then $D_{k-l}$ has one positive root and no negative root. If $k-l$ is even, then $D_{k-l}(x)$ has one positive root and one negative root.
\end{lemma}
\begin{proof}
Suppose $k-l$ is odd;  if $k$ is even then $l$ is odd. From
(\ref{equationDk-l}), with $-x\leftarrow x$, the numerator of
$D_{k-l}(x)$ can be written
\begin{equation}\label{A}
\sum_{t=0}^{l+1}(-1)^{t}a_t(-x)^{k+1-t}+(-1)^{l}l!=-\sum_{t=0}^{l+1}a_{t}x^{k+1-t}-l!.
\end{equation}
If $k$ is odd, then $l$ is even, and
\begin{equation}\label{B}
\sum_{t=0}^{l+1}(-1)^{t}a_t(-x)^{k+1-t}+(-1)^{l}l!=\sum_{t=0}^{l+1}a_{t}x^{k+1-t}+l!.
\end{equation}

\noindent By inspection of (\ref{A}), (\ref{B}) and employing
Descartes' rule,  $D_{k-l}(x)$ has no negative roots. Suppose $k-l$
is even; if $k$ is even then $l$ is even, and,
\begin{equation*}
\sum_{t=0}^{l+1}(-1)^{t}a_t(-x)^{k+1-t}+(-1)^{l}l!=-\sum_{t=0}^{l+1}a_{t}x^{k+1-t}+l!.
\end{equation*}

\noindent If $k$ is odd  then $l$ is odd, and
\begin{equation*}
\sum_{t=0}^{l+1}(-1)^{t}a_t(-x)^{k+1-t}+(-1)^{l}l!=\sum_{t=0}^{l+1}a_{t}x^{k+1-t}-l!.
\end{equation*}

\noindent By similar argument  $D_{k-l}(x)$ has one negative root.
Taking the $l$-th derivative of (\ref{equationFk}), it is easy to
see by Descartes' rule  that $D_{k-l}(x)$  has exactly one positive
root.
\end{proof}


\noindent Denote by $u_k$ the positive root of $D_k(x);$ for $k$
even, denote by $v_k$ the negative root of $D_k(x).$

\begin{theorem}We have the following results for the set $U$ and fixed
$l$:
\begin{enumerate}
\item Let $j=k-l.$  Then
\begin{equation*}
\lim_{j\to \infty}u_j=2.
\end{equation*}
 All of  the other complex roots of $D_j(x)$
are inside of $|z|<u_j.$ For $j$ even, we have
\begin{equation*}
\lim_{j\to \infty}v_j=-1
\end{equation*}
\item
If $j$ is odd, then $D_{j}(x)$ has one positive root and no
negative root. If $j$ is even, then $D_{j}(x)$ has one positive
root and one negative root.
\item
For $j\ge 2,$ we have $u_{j+1}>u_j.$
\item
There exists a even number $N_0,$ such that for even $n>N_{0},$ we have $v_{n+2}<v_{n}.$
\end{enumerate}
\end{theorem}
\begin{proof}
This theorem is proved by the following lemmas \ref{droot}--\ref{dnmon}.
\end{proof}
\begin{remark} The corresponding theorem has been proved in \cite{george2} for the first derivative and second derivative cases, .
\end{remark}
\begin{lemma}\label{dlim} Let $j=k-l,$ \ fixed $l.$ Then the positive
roots $u_j$ satisfy
\begin{equation*}
\lim_{j\to \infty}u_j=2.
\end{equation*}
All of the   the other complex roots of $D_j(x)$ are inside of open
disk $|z|<u_j.$ For $j$ even, the negative roots $v_j$ satisfy,
\begin{equation*}
\lim_{j\to \infty}v_j=-1.
\end{equation*}
\end{lemma}
\begin{proof} We have from (\ref{equationFk})
\begin{equation}\label{equationD}
D_{k-l}(x)-D_{k-1-l}(x)=(k-1)\cdots (k-l+1)x^{k-l-1}\left((x-2)k+2l\right).
\end{equation}
It follows that for any $a,$  $1<a<2,$
\begin{equation*}
\lim_{k\to \infty}D_{k-l}(a)-D_{k-1-l}(a)=-\infty.
\end{equation*}
Hence for any $a, \, 1<a<2,$ we have
\begin{equation*}
\lim_{k\to \infty}D_{k-l}(a)=-\infty.
\end{equation*}
It is easy to see from (\ref{equationD}) that
\begin{equation*}
\lim_{k\to \infty}D_{k-l}(2)=\infty.
\end{equation*}
Hence by the intermediate value theorem, $1< u_j<2$ for all  $j\ge j_0$  for sufficiently large $j_0.$
 \[\lim_{j\to \infty}u_j=2.\]

\noindent For $j$ even, we have from (\ref{equationFk})
\begin{equation}\label{Dkeven}
D_{k-l}(x)-D_{k-l-2}(x)=(k-2)\cdots (k-l+1)x^{k-l-2}h_{k}(x),
\end{equation}
where
\begin{equation}\label{Dh}
h_k(x)=(x^2-x-2)k^2+(-x^2+(l+1)x+2(2l+1))k-lx-2l(l+1).
\end{equation}
Hence if $a\le -1,$ we have from (\ref{Dkeven}), (\ref{Dh}),
\begin{equation*}
\lim_{k\to \infty}\left(D_{k-l}(a)-D_{k-l-2}(a)\right)=\infty.
\end{equation*}
For sufficiently large $k,$ if $-1<a<0,$  we have
\begin{equation*}
D_{k-l}(a)-D_{k-l-2}(a)<0.
\end{equation*}
Hence
for $j$ even, we have
\begin{equation*}
\lim_{j\to \infty}v_j=-1.
\end{equation*}

\noindent Notice for $j=k-l>0,$ we have
\begin{equation}\label{EquationDj}
D_j(x)=k(k-1)\cdots (k-l+1) x^{k-l}-\sum_{s=l}^{k-1}s(s-1)\cdots (s-l-1)x^{s-l}.
\end{equation}
Let $x_0=\rho e^{i\theta}$ be a complex zero of $D_j(x).$  By applying triangle inequality to (\ref{EquationDj}), we get $D_j(\rho)\le 0.$  We know that $D_j(x)<0$ if $0\le x\le u_j$ and $D_j(x)>0$ if $x>u_j.$  Since $\rho>0,$  we get $0<\rho <u_j.$
\end{proof}
\begin{lemma}\label{dpmon}
For $k\ge 2,$ we have $u_{k+1}>u_k.$
\end{lemma}
\begin{proof}
Solving
\begin{equation*}
D_k(x)-D_{k-1}(x)=l!\binom{k+l}{l}x^k-2l!\binom{k+l-1}{l}x^{k-1}=0,
\end{equation*}
we get
\begin{equation*}
x_k=\frac{2k}{k+l}=2-\frac{2l}{k+l}
\end{equation*}
Hence $x_k$ converges monotonically to $2$. We calculate
\begin{eqnarray}
 D_2(x_2)
&=& l!\left(\frac{(l+2)(l+1)}{2}\frac{4^2}{(l+2)^2}-(l+1)\frac{4}{l+2)}-1\right)\\
&=&l!\left(\frac{4(l+1}{l+2}-1\right)\nonumber\\
&=&l!\frac{3l+2}{l+2}\nonumber\\
&>& 0.\nonumber
\end{eqnarray}
Since $x_3>x_2,$  we obtain
\begin{equation}
D_3(x_3)=D_2(x_3)>D_2(x_2)>0.
\end{equation}
Hence $u_3>u_2.$  Inductively, we get $u_{k+1}>u_{k}.$
\end{proof}
\begin{lemma}\label{dnmon}
There exists an even number $N_0,$ such that for even $n>N_{0},$ we have $v_{n+2}<v_{n}.$
\end{lemma}
\begin{proof}
Solving
\begin{equation}\label{Dkn}
D_k(x)-D_{k-2}(x)=(k+l-2)\cdots (k+1)g_k(x)=0,
\end{equation}
where
\[
g_k(x)=(k+l)(k+l-1)x^2-k(k+l-1)x-2k(k-1),
\]
we get the negative root of (\ref{Dkn})
\begin{equation}\label{xk1}
x_{k}=\frac{k}{2(k+l)}-\frac{1}{2}\sqrt{\left(\frac{k}{k+l}\right)^2+\frac{8k(k-1)}{(k+l)(k+l-1)}}
\end{equation}
Consider the following function derived from (\ref{xk1})
\begin{equation}
f(x)=\frac{1}{2}(1-lx)-\frac{1}{2}\sqrt{(1-lx)^2+\frac{8(1-lx)\left(1-(l+1)x\right)}{(1-x)}}
\end{equation}
we find that
\begin{equation*}
f'(0)=\frac{5l}{2}>0.
\end{equation*}
Hence $f(x)$ is increasing on a neighborhood $V$ of $0.$

Since \begin{equation*}
\frac{1}{k+l}>\frac{1}{k+l+2},
\end{equation*}
we get
\begin{equation*}
x_{k+2}=f\left(1/(k+l+2)\right)<f\left(1/(k+l)\right)=x_{k}.
\end{equation*}

\noindent First we claim that there exists a sufficiently large even
number $k_{0},$  such that $v_{k_{0}}<v_{{k_0}-2}.$ Otherwise,
suppose there exists a $j_0,$  such that for all even number
$j>j_{0},$ $v_{j+2}\ge v_{j}.$  Since $D_k(-1) \rightarrow \infty$
as $k \rightarrow \infty,$ this contradicts the fact $\lim_{j\to
\infty}v_j=-1.$ Hence there exists an  even number $k_{0},$ such
that $v_{k_{0}}<v_{k_{0}-2}.$ It follows that
\begin{equation*}
D_{k_{0}}(x_{k_{0}})>0.
\end{equation*}
Otherwise, we have $v_{k_{0}}>v_{k_{0}-2},$  a contradiction. Since
$x_{k_{0}+2}<x_{k_{0}},$  we get
$D_{k_{0}+2}(x_{k_{0}+2})=D_{k_{0}}(x_{k_{0}})>0.$  It follows that
$v_{k_{0}+2}<v_{k_{0}}.$ Notice $\{x_{k}\}$ decreases to $-1$ also.
Inductively, we have that $v_{k+2}<v_{k}$ for $k$ sufficiently large
and even.
\end{proof}

\subsection{$V$ or integral case}\label{V}

\subsubsection {First Integral Case}
Now we consider the infinite sequence of polynomials $\{\int
F_k(x)\}$ of the first integral of the sequence $\{F_k(x)\}.$
\begin{definition}
We specify the following degree $j+1$ polynomial $I_j(x)$ to correspond with the first integral of $F_{j}(x).$

\begin{equation}\label{Ikx}
I_j(x)=\int F_j(x)=\frac{x^{j+1}}{j+1}-\frac{x^j}{j}-\cdots-x-1,
\end{equation}
for all $j\ge 1.$
\end{definition}
\begin{theorem} The roots of $I_k(x)$ satisfy the following
properties,
\begin{enumerate}
\item $I_k(x)$ has a positive simple  root $\phi_k$ satisfying $2<\phi_k<3$.
\item  For $k\ge 2,$ we have $\phi_{k+1}<\phi_k.$
\item
\begin{equation*}
\lim_{j\to \infty}\phi_k=2.
\end{equation*}
\item
If $k$ is odd, then
\begin{enumerate}
\item  $I_k(x)$ has a negative simple  root $\theta_k$ satisfying $-2<\theta_k<-1$.
\item
\begin{equation*}
\lim_{j\to \infty}\theta_k=-1.
\end{equation*}
\item
$\theta_k>\theta_{k-2}$ for $k\ge 17.$
\end{enumerate}
\item For $k$ even, $I_k(x)$ has no negative root.
\end{enumerate}
\end{theorem}
\begin{proof}
We prove this theorem in the following lemmas \ref{uni}--\ref{enn}.
\end{proof}
\begin{lemma}\label{uni} $I_k(x)$ has a positive  simple  root $\phi_k$ satisfying $2<\phi_k<3$.
If $k$ is odd, then $I_k(x)$ has a negative simple root $\theta_k$
satisfies $-2<\theta_k<-1$.
\end{lemma}
\begin{proof} From Descartes' Rule, we get that the number of possible positive roots for each $I_k(x)$ is $1$.
If $a=2,$
\begin{equation*}
I_1(2)=\frac{2^2}{2}-2-1=-1<0
\end{equation*}
We find that for $k>1,$
\begin{equation*}
I_k(2)-I_{k-1}(2)=\frac{-2^{k+1}}{k(k+1)}<0
\end{equation*}
Then for all $k\ge 1,$ we have $I_k(2)<0.$  Hence the positive root
$\phi_k>2.$ If $a=3,$ then
\begin{equation*}
I_1(3)=\frac{3^2}{2}-3-1=\frac{1}{2}>0
\end{equation*}
We have that for $k\ge 2,$
\begin{equation}
I_k(3)-I_{k-1}(3)=\frac{3^{k}(k-2)}{k(k+1)}\ge 0
\end{equation}
Then for all $k\ge 1,$ we have $I_k(3)>0.$  Hence the positive root
$\phi_k$ satisfies $2<\phi_k<3.$ If $I'_k(\phi_k)=0,$ then by
\cite{miller}, $1< \phi_k<2.$  Hence $I'_k(\phi_k)\ne 0.$ Therefore,
$\phi_k$ is a simple root of $I_k(x).$ For $k$ odd, we can get that
the number of variation for signs $I_k(-x)$ is $k.$  Then by
Descartes' Rule, we know the possible  number of negative roots for
$I_k(x)$ is $k, \, k-2, \dots, \, k-2t,\dots, 1.$  By \cite{miller},
we know $I'_k(x)$ only has one real root $b_k$ for $k$ odd.  It
follows that $I_k(x)$ is increasing if $x>b_k$ and decreasing if
$x<b_k.$  Hence we get the number of negative real roots for
$I_k(x)$ is $1.$

\noindent If $k=1$ then
\begin{equation*}
I_1(-1)=-1+1+\frac{1}{2}=\frac{1}{2}>0.
\end{equation*}
If $k=3$  then
\begin{equation*}
I_3(-1)=-1+1-\frac{1}{2}+\frac{1}{3}+\frac{1}{4}>0.
\end{equation*}
If $k=5$ then
\begin{equation}\label{I5}
I_5(-1)=-1+1-\frac{1}{2}+\frac{1}{3}-\frac{1}{4}+\frac{1}{5}+\frac{1}{6}=-\frac{1}{20}<0.
\end{equation}
If $k>5$ and $k$ is odd then
\begin{equation}\label{Ik1}
I_k(-1)=\frac{1}{k+1}+\sum_{l=1}^k\frac{(-1)^{l-1}}{l}-1\le -1+\ln(2)+\frac{2}{k+1}<0
\end{equation}
For $k$ odd, from \cite{miller}, $I'_k(x)<0$ for $x<0,$ so $I_k(x)$
is decreasing for $x<0.$ Hence from (\ref{I5})and (\ref{Ik1}) for
all $k\ge 5$ and $k$ odd, the negative real root $\theta_k$ of
$I_k(x).$ satisfying $\theta_k<-1.$

\noindent Next we show $-2<\theta_k.$ From (\ref{Ikx}), we obtain
\begin{equation*}
I_k(x)-I_{k-2}(x)=\frac{x^{k+1}}{k+1}-\frac{x^{k}}{k}-2\frac{x^{k-1}}{k-1}.
\end{equation*}
Solving
\begin{equation}
\frac{x^2}{k+1}-\frac{x}{k}-\frac{2}{k-1}=0,
\end{equation}
yields the negative root,
\begin{equation}\label{xk1xk2}
x_{k1}=\frac{1}{k}-\sqrt{\left(\frac{1}{k}\right)^2+\frac{8}{(k+1)(k-1)}}\frac{k+1}{2},
\end{equation}
It can be shown by direct calculation  that for $k$ odd and $k\ge 7,$
$-1>x_{k1}>-2.$

\noindent That implies that for $k$ odd and $k\ge 7,$
\begin{equation}
I_k(-2)-I_{k-2}(-2)>0
\end{equation}
We know
\begin{equation*}
I_5(-2)=\frac{221}{15}>0.
\end{equation*}
Hence for all $k\ge 5$ and $k$ odd, we have $I_k(-2)>0.$
Therefore we get $-2<\theta_k<-1.$

\noindent If $I'_k(\theta_k)=0,$ then by \cite{miller},
$-1<\theta_k<0.$ Hence $I'_k(\theta_k)\ne 0.$  It follows that
$\theta_k$ is a simple root of $I_k(x).$
\end{proof}
\begin{lemma}\label{mon}
Let $\phi_k$ be the positive root of $I_k(x).$  Then for $k\ge 2,$ we have $\phi_{k+1}<\phi_k.$
\end{lemma}
\begin{proof} Denoted by $b_k$ the positive real root of $I'_k(x).$  By \cite{miller}, we know $1<b_k<2.$  Hence $I_i(x),$ $i\ge 2,$  is increasing if $x>2.$  It's easy to see
that $I_{i}>I_{i-1}$ if $x>2+\frac{2}{i}$ and $I_{i}<I_{i-1}$ if $x<2+\frac{2}{i}.$  Notice $2+\frac{2}{i}$ converges to $2$ decreasingly. From $I_3(2+2/3)<0,$  we get
$\phi_3<\phi_2.$  Suppose for all $2<i\le k,$  we have $I_{i}(2+2/i)<0.$  Then since
$I_{k-1}(2+2/k)=I_{k}(2+2/k)<0$ and $I_{k+1}$ is increasing if $x>2+2/(k+1),$  we get $I_{k+1}(2+2/(k+1))<0.$
 We know $I_{k+1}>I_{k}$ if $2+2/(k+1)<x<3.$ We get $I_k(\phi_{k+1})<0.$  Hence for $k\ge 2,$ $\phi_{k+1}<\phi_k.$
\end{proof}
\begin{lemma}\label{lim2}
\begin{equation*}
\lim_{k\to \infty}\phi_k=2.
\end{equation*}
\end{lemma}

\begin{proof}
For any $k,$
\begin{eqnarray}\label{IkIk-1}
I_k(x)-I_{k-1}(x)&=&x^k\left(\frac{x}{k+1}-\frac{2}{k}\right)\\
                 &=&\frac{x^k[(x-2)k-2]}{k(k+1)}.\nonumber
\end{eqnarray}
If $a>2,$ then for sufficiently large $k$,
\begin{equation}\label{(a-2)k}
(a-2)k-2>1.
\end{equation}
We know
\begin{equation}\label{xk}
\lim_{k\to \infty}\frac{x^k}{k(k+1)}=\infty.
\end{equation}
Hence employing  (\ref{(a-2)k}), (\ref{xk}) in (\ref{IkIk-1}), for any $a>2,$ yields
\begin{equation}\label{Ikinf}
\lim_{k\to \infty} I_k(a)-I_{k-1}(a)=\infty.
\end{equation}

\noindent Notice for any $k>2,$
\begin{equation}\label{Ikinf1}
I_k(x)=\sum_{l=3}^{k}\left(I_l(x)-I_{l-1}(x)\right)+I_2(x).
\end{equation}
It follows from (\ref{Ikinf}), (\ref{Ikinf1}) that for any $a>2,$
\begin{equation*}
\lim_{k\to \infty}I_k(a)=\infty.
\end{equation*}
If $a=2,$ we have

\begin{equation*}
I_k(2)-I_{k-1}(2)=\frac{-2.2^k}{k(k+1)}.
\end{equation*}
Then by a similar argument as above,
\begin{equation*}
\lim_{k \to \infty}I_k(2)=-\infty.
\end{equation*}
By the Mean Value Theorem, we obtain
\begin{equation*}
\lim_{k\to \infty}\phi_k=2.
\end{equation*}
\end{proof}

\begin{lemma}\label{limn1}  Let $k$ be a odd number and $\theta_k$ be the negative root of $I_k(x).$ Then
\begin{equation}
\lim_{k\to \infty}\theta_k=-1.
\end{equation}
Moreover, for $k>17$ and $k$ is odd, $\theta_k>\theta_{k-2}.$
\end{lemma}
\begin{proof}
For $a<-1,$  we have from (\ref{Ikx})

\begin{eqnarray}\label{Ika-Ik-2}
I_k(a)-I_{k-2}(a)&=&a^{k-1}\left(\frac{a^2}{k+1}-\frac{a}{k}-\frac{2}{k-1}\right)\\
                 &=&a^{k-1}\frac{k^2(a^2-a-2)-(a^2+2)k+a}{(k+1)k(k-1)}\nonumber
\end{eqnarray}
For $a<-1$ and $k$ sufficiently large, we have from (\ref{Ika-Ik-2})
\begin{equation*}
k^2(a^2-a-2)-(a^2+2)k+a>1
\end{equation*}
Since for $a<-1$ and $k$ odd, by similar argument as lemma \ref{lim2},
\begin{equation*}
\lim_{k\to \infty}\frac{a^{k-1}}{(k+1)k(k-1)}=\infty,
\end{equation*}
we get
\begin{equation*}
\lim_{k\to \infty} I_k(a)-I_{k-2}(a)=\infty.
\end{equation*}

\noindent Then by writing $I_k(x)$ as telescoping sum,
 $a<-1$, $k$ odd,  it follows that

\begin{equation}\label{Ika1}
\lim_{k\to \infty}I_k(a)=\infty.
\end{equation}



\noindent Substituting  $a=-1$ in (\ref{Ikx}) gives
\begin{equation*}
I_k(-1)=\frac{(-1)^{k+1}}{k+1}+H_k(-1)-1,
\end{equation*}
where $H_k(x)$ is the standard alternating sum.

\noindent Hence \begin{equation}\label{Ik-1} \lim_{k\to
\infty}I_k(-1)=\ln(2)-1<0.
\end{equation}
It follows that from Mean Value Theorem, (\ref{Ika1}), (\ref{Ik-1}),
 \begin{equation*}
\lim_{k\to \infty}\theta_k=-1.
\end{equation*}
A calculator check  with $k=17$ in (\ref{xk1xk2}) yields  \[I_k(x_{k1})=-0.0337812682<0.\]

\noindent From (\ref{xk1xk2}) we write
 \begin{equation*}
f(x)=-\sqrt{(1+x)^2+8\frac{1+x}{1-x}}+(1+x)
\end{equation*}
Taking the derivative of $f(x)$ gives
\begin{equation*}
f'(x)=-\frac{1}{2}\frac{2(1+x)+16/(1-x)^2}{\sqrt{(1+x)^2+8(1+x)/(1-x)}}+1
\end{equation*}
It's easy to check that for $0<x<1,$  $f'(x)<0.$  $f(x)$ is
decreasing for $0<x<1.$ Since  $1/(k+2)<1/k,$  we get for $k\ge 7,$
\begin{equation}
x_{k}=\frac{\frac{1}{k}-
\sqrt{\left(\frac{1}{k}\right)^2+\frac{8}{(k+1)(k-1)}}}{\frac{2}{k+1}}
<\frac{\frac{1}{k+2}-\sqrt{\left(\frac{1}{k+2}\right)^2+\frac{8}{(k+3)(k+1)}}}{\frac{2}{k+3}}=x_{k+2}.
\end{equation}

\noindent Hence $x_{k}$ increases to $-1.$  Denote by $\theta_k$ the
negative real root of $I_k(x).$ Since $I_{17}(x_{17})<0$,  we get
$\theta_{17}<x_{17}.$ It follows that $I_{19}(\theta_{17})>0$ since
$I_{19}(x)>I_{17}(x)$ when $x<x_{17}.$   Hence
$\theta_{19}>\theta_{17};$
 it follows that  $\theta_k>\theta_{k-2}$ for $k\ge 17.$
\end{proof}
\noindent It is noted that for  $1<a<2$,   using  similar methods,
we can get
\begin{equation*}
\lim_{k\to \infty}I_k(a)=-\infty.
\end{equation*}
\begin{lemma}\label{enn} For $k$ even, the integral $I_k(x)$, $(\ref{Ikx})$,  has no negative root.
\end{lemma}
\begin{proof}
Let $k=2l,$  $x=-a$ for $0<a<1.$
By rewriting (\ref{Ikx}) we get
 \begin{equation}
I_k(x)=-\frac{a^{2l+1}}{2l+1}-\frac{a^{2l}}{2l}+\sum_{t=2}^{l}a^{2(t-1)}\left(\frac{a}{2t-1}-\frac{1}{2t-2}\right)+a-1<0.
\end{equation}
Hence, for $k$ even, $I_k(x)$ has no negative root on $-1<x<0.$
It is easy to check that $I_k(-1)<0.$

\noindent By \cite{miller}, for $k$ even, $I'_k(x)$ has a negative
root $r_k$ satisfying $-1<r_k<0.$  Hence $I_k(x)$ is increasing on
$-\infty<x<-1$ so that    for $k$ even $I_k(x)<0.$ Therefore, for
$k$ even, $I_k(x)$ has no negative root.
\end{proof}
\begin{lemma} For any $k\ge 2$,  the complex zeros of $I_k(z)$   satisfy the inequality $|z|<\phi_k<3.$
\end{lemma}
\begin{proof} Let $z_{0}=re^{i\theta}$ be a complex root of $I_k(z).$  Using the triangle inequality we obtain
\begin{equation}
I_k(r)\le 0.
\end{equation}
Note that equality holds only at $\theta=0,$ i.e $z_{0}=\phi_k.$  Since $I_k(x)<0$  for $0<x<\phi_k<3$ and $x$ real, we get $r<\phi_k<3.$
\end{proof}
\begin{lemma} If $-1<a<1,$ then
\[|I_k(x)|\leq  \frac{1}{1-|x|}.\]
Moreover,
\begin{equation*}
\lim_{k\to \infty}I_k(x)=-1+\ln(1-x).
\end{equation*}
\end{lemma}
\begin{proof}
If  $ -1<a<1,$
then
\begin{eqnarray*}
|I_k(x)|&\leq &\sum_{l=0}^k\left|\frac{x^{l+1}}{l+1}\right|+1\\
         &\leq &\sum_{l=0}^k\left|{x^{l+1}}\right|+1\nonumber\\
         &\leq & \frac{1}{1-|x|}.\nonumber
\end{eqnarray*}
The Taylor series expansion for  $I_k(x)$ with $-1<x<1,$ yields
\begin{equation*}
\lim_{k\to \infty}I_k(x)=-1+\ln(1-x).
\end{equation*}
\end{proof}

\subsubsection{General Case}Now we consider the infinite sequence of polynomials $\left\{\overbrace{\int\int\cdots\int}^{l+2} F_{k}(x)\right\}$ of the $(l+2)$-th integral  of the  sequence  $\{F_k(x)\}.$
\begin{definition}
For $0<l<k,$  We specify the following degree $k+1$ polynomial $H_k(x)$ to correspond with the $(l+2)$-th integral of $F_{k-l-1}(x).$
\begin{eqnarray}\label{equationhkx}
H_k(x)&=&\overbrace{\int\int\cdots\int}^{l+2} F_{k-l-1}(x)\\ &=&\frac{x^{k+1}}{(l+2)!\binom{k+1}{l+2}}-\sum_{t=l+2}^k\frac{x^t}{(l+2)!\binom{t}{l+2}}-\sum_{s=0}^{l+1}\frac{x^s}{s!}\nonumber
\end{eqnarray}
Let $\alpha_k$ be the positive root of $H_k(x).$  For $k$ odd, denote by $\beta_k$ the negative real root of $H_k(x).$
\end{definition}
\noindent We have the following
\begin{theorem} The roots of  $H_k(x)$ satisfy the following
properties,
\begin{enumerate}
\item
\begin{equation*}
\lim_{k\to \infty}\alpha_k=2.
\end{equation*}
Except $\alpha_k,$  all the other complex roots are inside $\{z:|z|<\alpha_k\}.$
For $k$ odd, we have
\begin{equation*}
\lim_{k\to \infty}\beta_k=-1.
\end{equation*}

\item
For  sufficiently large even $k,$ for any $x<0,\, H_k(x)<0,$ i.e $H_k(x)$ has no negative real roots.

\item
For  sufficiently large odd $k,$  for any $x<0,$
$H_k(x)$ has one negative root.

\item
 $\alpha_{j+1}<\alpha_j,$  $\forall j\ge l+3,$

\item
there exists odd $N_{0},$  such that for all odd $n\ge N_{0},$  we have $\beta_{n+2}>\beta_n.$
\end{enumerate}
\end{theorem}
\begin{proof}
The theorem is proved using lemmas~\ref{posr}--\ref{nmon}.
\end{proof}
\begin{lemma}\label{posr}
\begin{equation*}
\lim_{k\to \infty}\alpha_k=2.
\end{equation*}
Except $\alpha_k,$  the other complex roots are inside $\{z:|z|<\alpha_k\}.$
For $k$ odd, we have
\begin{equation*}
\lim_{k\to \infty}\beta_k=-1.
\end{equation*}
\end{lemma}
\begin{proof}
The proof uses similar idea as the previous section with some
differences, we include for completeness.

\begin{equation}\label{eqhkhk-1}
H_k(x)-H_{k-1}(x)
=\frac{x^k}{(l+3)!\binom{k+1}{l+3}}\left((x-2)k-lx-x-2)\right).
\end{equation}

\noindent It follows that for $a>2,$
\begin{equation*}
\lim_{k\to \infty}\left(H_k(a)-H_{k-1}(a)\right)=\infty.
\end{equation*}
Hence for $a>2,$
\begin{equation*}
\lim_{k\to \infty}H_k(a)=\infty.
\end{equation*}
It's easy to prove that
\begin{equation*}
\lim_{k\to \infty}H_k(2)=-\infty.
\end{equation*}
Hence,
\begin{equation*}
\lim_{k\to \infty}\alpha_k=2.
\end{equation*}
Let $z=re^{i\theta}.$  Then by triangle inequality,
\begin{equation}\label{hk0}
H_k(r)\le 0
\end{equation}
Equality in (\ref{hk0}) holds only at $\theta=0;$
it follows that $r<\alpha_k.$  Since $z=0$ is not the root of $H_k(z),$  we have $0<r<\alpha_k.$

\noindent If $k$ is odd, then
\begin{equation}\label{hkhk-2}
H_k(x)-H_{k-2}(x)
=\frac{x^{k-1}}{(l+4)!\binom{k+1}{l+4}}h_k(x),
\end{equation}
where
\begin{equation}\label{hkx}
h_k(x)=(x^2-x-2)k^2-\left((2l+3)x^2+(l+1)x+2\right)k+\left((l+1)(l+2)x^2+(l+2)x\right).
\end{equation}
Hence, if $a<-1$ and $k$ odd, employing (\ref{hkhk-2}), (\ref{hkx})  we have
\begin{equation}\label{hkahka-2}
\lim_{k\to \infty}\left(H_k(a)-H_{k-2}(a)\right)=\infty.
\end{equation}
It follows from (\ref{hkahka-2}) that
\begin{equation}\label{hkinfty}
\lim_{k\to \infty}H_k(a)=\infty.
\end{equation}
For $k$ odd, it is easy to see from (\ref{equationhkx}) that for sufficiently large $k,$
\begin{equation}\label{hk-1}
H_k(-1)<0
\end{equation}
Denote by $\beta_k$ the negative real root of $H_k(x).$
 We have from (\ref{hkinfty}), (\ref{hk-1})
 \begin{equation*}
\lim_{k\to \infty}\beta_k=-1.
\end{equation*}

\end{proof}
\begin{lemma}\label{negr}
For  sufficiently large even $k,$  for any $x<0,$ $H_k(x)<0.$

\end{lemma}
\begin{proof}
 This result was shown for the first integral ($l=-1$ in (\ref{equationhkx})) in lemma \ref{enn}.
 Now we consider the case $l\ge 0$ in (\ref{equationhkx}).

\noindent For $k$ and $l$ both even, we obtain
\begin{equation*}
H'_{k}(-1)=A+ B+C-1,
\end{equation*}

\noindent where \begin{equation*} A=\frac{1}{(l+1)!\binom{k}{l+1}}+
\frac{1}{(l+1)!\binom{k-1}{l+1}},
\end{equation*}
\begin{equation*}
B=\sum_{d=l+3}^{k-2}\frac{l+1}{(l+2)!\binom{d}{l+2}},
\end{equation*}
\begin{equation*}
C=\sum_{s=1}^{l/2}\frac{2s-1}{(2s)!}.
\end{equation*}

\noindent  We note that $H'_k(-1)<0.$ Since$A>0,$ $B>0,$ and $C\ge
\frac{1}{2},$ this implies for $k$  sufficiently large even $k$ and
$l$ even,
\begin{equation}\label{equationhk-1}
|H'_{k}(-1)|<\frac{1}{2}.
\end{equation}
The same result (\ref{equationhk-1}) holds with a similar proof   in the case of  odd $l$ and for  sufficiently large even $k.$
Let $\theta_k$ be the negative root of $H'_k(x)$ and let $\gamma_k$ be the negative root of $H_{k}^{(3)}(x).$ We know from lemma~\ref{posr},
\begin{equation}
\lim_{k\to \infty}\theta_k=-1,\, \lim_{k\to \infty}\gamma_k=-1.
\end{equation}
 Notice \[H_{k}(\theta_k)=\int_{0}^{\theta_k} H'_{k}(x)dx-1\]
 and $H'_k(x)$ is decreasing on $x<0$ so $|H'_k(x)|<1$ since $H'_k(0)=-1.$  $H'_k(x)$ is concave down on $\gamma_k<x<0$ since $H_{k}^{(3)}(x)<0$ on $\gamma_k<x<0.$

\noindent Hence for  sufficiently large even $k,$ if
$\theta_k<\gamma_k,$ we obtain
\begin{eqnarray}\label{hktheta}
H_{k}({\theta_k})&=&\int_{\gamma_k}^{\theta_k}H'_{k}(x)dx+\int_{0}^{\gamma_k}H'_k(x)dx-1\\
&<&|\theta_k-\gamma_k|+\frac{1}{2}\left(|H'_k(\gamma_k)|+1\right)-1.\nonumber
\end{eqnarray}

\noindent For $l=0,$   we know that $-2<\theta_k<-1$  and
$-1<\gamma_k<0$.
Write $\gamma_k=-a_k$  and write $k=2t.$  Then by taking the derivative of (\ref{equationhkx}),
\begin{equation}\label{h'kgammak}
H'_k(\gamma_k)=-1+\sum_{s=1}^{t-1}{a_k}^{2s-1}\left(\frac{1}{2s-1}-\frac{a_k}{2s}\right)+\frac{a_{k}^{2t-1}}{2t-1}+\frac{a_k^{2t}}{2t}.
\end{equation}
Hence since $-1<H'_k(\gamma_k)<0$ and by inspection of (\ref{h'kgammak})
\begin{equation}\label{h'k1}
|H'_k(\gamma_k)|<\left|-1+a_k-\frac{{a_k}^2}{2}\right|<1.
\end{equation}
It follows for  sufficiently large $k$ from (\ref{hktheta}),(\ref{h'k1})
\[H_{k}(\theta_k)<0.\]
Therefore $H_k(x)<0$ for all $x<0.$

\noindent For $l>0,$  we have $\gamma_k<-1.$  If
$\theta_k<\gamma_k,$ then $|H'_k(\gamma_k)|<|H'_k(-1)|<\frac{1}{2}$
in (\ref{equationhk-1}). Hence for sufficiently large even $k$, we
get
\begin{equation}
H_{k}(x)<0.
\end{equation}

\noindent If $\gamma_k<\theta_k,$  then
\begin{equation}
H_{k}(\theta_k)<\frac{1}{2}|\theta_k|-1<0.
\end{equation}
It follows for  sufficiently large even $k$, we get
\begin{equation*}
H_{k}(x)<0.
\end{equation*}

\end{proof}
\begin{lemma}\label{nor}
For sufficiently large odd $k,$ for any $x<0,$
$H_k(x)$ has exactly one negative root.
\end{lemma}
\begin{proof}
By Lemma \ref{negr}, we know $H'_k(x)<0$ for $x<0.$  Hence $H_k(x)$ is decreasing on $(-\infty, 0).$  Since
$H_k(0)=-1$ and $\lim_{x\to -\infty}H_k(x)=\infty,$  we get that $H_k(x)$ has only one root on $( -\infty, 0).$
\end{proof}
\noindent Now we study the monotonicity of the positive root
$\alpha_k$ of $H_k(x)$ in the following
\begin{lemma}\label{pmon}
For all $j\ge l+3,$ where $l\ge -1$ is a fixed integer, we have $\alpha_{j+1}<\alpha_j.$
\end{lemma}
\begin{proof}
Solving for the zero of (\ref{eqhkhk-1}) for $k=l+3$  yields
the intersection point $x=l+4=k+1.$ Next we show $H_k(x)<0$ at the
intersection point $x=k+1,$

\begin{eqnarray}
 H_k(k+1)
&=&\frac{(k+1)^{k+1}}{(k-1)!\binom{k+1}{k-1}}-\sum_{t=1}^{k}\frac{(k+1)^t}{t!}-1\\
&=& \frac{(k+1)^k}{k!}-\sum_{t=1}^{k-1}\frac{(k+1)^t}{t!}-1\nonumber\\
&=& \frac{(k+1)^{k-1}}{(k-1)!}\left(\frac{k+1}{k}-1\right)-\sum_{t=1}^{k-2}\frac{(k+1)^t}{t!}-1 \nonumber\\
&=& \frac{(k+1)^{k-1}}{k!}-\sum_{t=1}^{k-2}\frac{(k+1)^t}{t!}-1\nonumber\\
&=&\frac{(k+1)^{k-2}}{(k-2)!}\left(\frac{-k^2+2k+1}{k(k-1)}\right)-\sum_{t=1}^{k-3}\frac{(k+1)^t}{t!}-1\nonumber\\
&<& 0.\nonumber
\end{eqnarray}
The lemma follows the similar argument as  lemma (\ref{dpmon}).
\end{proof}
\noindent We now consider the monotonicity of the negative root
$\beta_k$ of $H_k(x)$ in the following
\begin{lemma}\label{nmon}
There exists odd $N_{0},$  such that for all odd $n\ge N_{0},$  we have $\beta_{n+2}>\beta_n.$
\end{lemma}
\begin{proof}


\noindent Solving the zero of (\ref{hkhk-2})
we get the negative real root
\begin{equation}\label{xkHk}
x_{k}=\frac{1}{2}\left({\frac{k+1}{k-l-1}-\sqrt{\frac{(k+1)^2}{(k-l-1)^2}+\frac{8(k+1)k}{(k-l-1)(k-l-2)}}}\right).
\end{equation}

\noindent We consider the function derived from (\ref{xkHk})
\begin{equation}
f(x)=1+(l+2)x-\sqrt{ \left(1+(l+2)x\right)^2+8\frac{\left(1+(l+2)x\right)\left(1+(l+1)x\right)}{1-x}}.
\end{equation}

\noindent Taking the derivative of $f(x)$ gives
\begin{equation}\label{eqf'x}
f'(x)=(l+2)-\frac{A+B+C}{2}\left(\sqrt{\left(1+(l+2)x\right)^2+8\frac{\left(1+(l+2)x\right)\left(1+(l+1)x\right)}{1-x}}\right)^{-1},
\end{equation}
where
\begin{equation*}
A=2\left(1+(l+2)x\right)(l+2),
\end{equation*}
\begin{equation*}
B=8(l+2)\frac{1+(l+1)x}{1-x},
\end{equation*}
\begin{equation*}
C=8\left(1+(l+2)x\right)\frac{(l+1)(1-x)+\left(1+(l+1)x\right)}{(1-x)^2}.
\end{equation*}
Substituting $x=0$ in (\ref{eqf'x}) gives
 \begin{equation}\label{f'0}
f'(0)=-2(l+2)<0.
\end{equation}
It follows from (\ref{f'0}) that there exists a neighborhood $V$ of $0,$  such that $f(x)$ is decreasing on $V.$

\noindent Since \begin{equation} \frac{1}{k-l-1}>\frac{1}{k+2-l-1},
\end{equation}
we have
\begin{equation}
x_{k}=\frac{1}{2}f\left({1}/{(k-l-1)}\right)<\frac{1}{2}f\left({1}/{(k+2-l-1)}\right)=x_{k+2}.
\end{equation}
It's easy to see that
\begin{equation}
\lim_{k\to \infty}x_{k}=-1.
\end{equation}

\noindent We claim that there exists a sufficiently large  odd
number $j_{0},$ such that $\beta_{j_{0}+2}>\beta_{j_{0}}.$
Otherwise, suppose there exists a $k_0,$  such that for all odd
number $n>k_{0},$  we always have $\beta_{n+2}\le \beta_{n}.$  This
contradicts the fact $\lim_{k\to \infty}\beta_k=-1.$ It follows that
\begin{equation}
H_{j_{0}}(x_{j_{0}})<0.
\end{equation}
Otherwise suppose $H_{j_{0}}(x_{j_{0}})>0.$ Since $H_{j}(x)$ is
decreasing on $x<0,$  we get $\beta_{j_{0}}>x_{j_{0}}.$  Since
$H_{j_{0}+2}(\beta_{j_{0}})<0,$  we get
$\beta_{j_{0}+2}<\beta_{j_{0}},$  a contradiction. Then the lemma
follows the similar arguments as lemma~\ref{limn1}.
\end{proof}


\begin{thebibliography}{10}
\bibitem{fd1} Franc\c{o}is Dubeau, {\em On $r$-Generalized
Fibonacci Numbers,} The  Fibonacci Quarterly, \textbf{27.3}, (1989),
pp.221-229. \vspace{5mm}

\bibitem{fd2}  \bysame, {\em The Rabbit Problem
Revisited,} The  Fibonacci Quarterly, \textbf{31.3}, (1993),
pp.268-273. \vspace{5mm}
\bibitem{flores} Ivan Flores, {\em Direct Calculation of $k$-Generalized Fibonacci Numbers},
 The Fibonacci Quarterly, \textbf{5.3}, (1967), pp. 259-266.
 \vspace{5mm}

 \bibitem{george3} George Grossman,
 {\em Recurrence relations and combinatorial identities} pre-print, 2006. \vspace{5mm}

 \bibitem{george4}\bysame,
 {\em Polynomial representation of binomial coefficients},
  pre-print, 2005.
 \vspace{5mm}

 \bibitem{george5} \bysame,  {\em Linear recurrence relations and the binomial coefficients},
 in Proceedings of XII$^{th}$ CZECH-POLISH-SLOVAK
Mathematical School by the Faculty of Education of University J. E.
Purkyn$\breve{e}$, \'{U}st\'{i} nad Labem, $Hublo\breve{s}$, June
2-4, 2005, pp. 111-119. \vspace{5mm}

\bibitem{sumid1}\bysame, Akalu Tefera and Aklulu Zeleke,
{\em Summation Identities for Representation of Certain Real
Numbers}, International Journal of Mathematics and Mathematical
Sciences(e-journal), Volume 2006 , Article ID 78739, 8 pages.
 \vspace{5mm}

\bibitem{fib2}\bysame, Akalu Tefera and Aklulu Zeleke, {\em On proofs of certain
combinatorial identities}, pre-print. \vspace{5mm}

\bibitem{sumid2}\bysame and Aklilu
Zeleke, {\em On linear recurrence relations and combinatorial
identities},  Journal of Concrete and Applicable Mathematics, Vol.
1, (2003), No. 3, pp. 229-245, Nova Science Publishers. \vspace{5mm}

\bibitem{george1}\bysame,
\newblock{\em{Fractal construction by orthogonal projection using the Fibonacci
sequence}},
\newblock {The Fibonacci Quarterly} {\bf35}, (1997), no. 3, 206-224.
\vspace{5mm}

\bibitem{george2}
\bysame and Sivaram Narayan,
\newblock{\em{On the characteristic polynomial of the $jth$ order Fibonacci
sequence}},
\newblock{ Applications of Fibonacci numbers, Vol. 8}, (1999), pp. 165-177
\newblock{Kluwer Acad. Publ.}, Dordrecht.
\vspace{5mm}

\bibitem{miles}
E. P. Miles Jr.,
\newblock {\em On Generalized Fibonacci Numbers and Associated Matrices},
\newblock { The Amer. Math. Monthly} {\bf67}, (1960), no. 8, pp.
745-752. \vspace{5mm}

\bibitem{miller}
M. Miller,
\newblock {\em On Generalized Fibonacci Numbers.},
\newblock { The Amer. Math. Monthly} {\bf78}, (1971), no. 10, pp.
1108-1109. \vspace{5mm}

\bibitem{truss} Alfred S. Posamentier and Ingmar Lehman, Afterword by
Herbert A. Hauptman, Nobel laureate. {\em The (Fabulous)Fibonacci
Numbers,} Prometheus Books, 2007. \vspace{5mm}

\bibitem{polszeg} G. P\'{o}lya and G. Szeg\"{o}, {\em Aufgaben \ and \
Lehrs\"{a}tze,} New York, Dover Publications, 1954.
\vspace{5mm}
\end{thebibliography}
\end{document}